\documentstyle[12pt]{article}

  \newtheorem{thm}{Theorem}[section]
  \newtheorem{rk}[thm]{Remark}
  \newtheorem{prop}[thm]{Proposition}
  \newtheorem{clly}[thm]{Corollary}
  \newtheorem{lemma}[thm]{Lemma}
  \newtheorem{defi}[thm]{Definition}
  \newtheorem{nota}[thm]{Notation}
  \newtheorem{exam}[thm]{Example}
  \newcommand{\pf}{\vspace{.2in}\hspace{.6cm}{\bf Proof$\,$:\ }}
  \newcommand{\bea}[1]{\begin{eqnarray}\label{#1}}
  \newcommand{\eea}{\end{eqnarray}}
  \newcommand{\beas}{\begin{eqnarray*}}
  \newcommand{\eeas}{\end{eqnarray*}}
  \newcommand{\lrangle}{\rangle_{{\scriptscriptstyle{L}}}}
  \newcommand{\rrangle}{\rangle_{{\scriptscriptstyle{R}}}}
  \newcommand{\cp}{\hbox to 1.8ex{$\times \kern-.45ex\vrule height1.1ex
depth0pt width0.45truept$\hfill}}

  \newcommand{\TT}{{\bf{T}}}
  \newcommand{\ZZ}{\mbox{$Z\!\!\!Z$}}
  
  \newcommand{\RR}{\mbox{$I\!\!R$}}
  \newcommand{\CC}{\mbox{$I\!\!\!C$}}

  \newcommand{\qed}{\hfill {\bf Q.E.D.} \vspace*{.1in}}
  \newcommand{\comp}{{\scriptstyle \circ}}
  \newcommand{\aut}{\mbox{Aut}}
  \newcommand{\proj}{\mbox{Proj}}
  \newcommand{\cpic}{\mbox{CPic}}
  \newcommand{\pic}{\mbox{Pic}}
  \newcommand{\gin}{\mbox{Gin}}

  \newcommand{\cstar}{\mbox{$C^*$}}
  
  \newcommand{\comsp}{\mbox{,\hspace{.2in}}}

  \newcommand{\dc}{D_{\mu\nu}^{c}}
  
  \newcommand{\mca}{M^c_{\alpha_{\mu\nu}}}
  \newcommand{\mc}{M^{c}}

  \newcommand{\hcb}{Hilbert \cstar-bimodule }
  \newcommand{\hcbs}{Hilbert \cstar-bimodules }
  \newcommand{\ct}{\mbox{$C({\bf{T}}^2)$}}
  \newcommand{\al}{\alpha_{\mu\nu}}
  \newcommand{\cpcong}{\stackrel{cp}{\cong}}

  \title{Hilbert \cstar-bimodules over commutative \cstar-algebras and an
isomorphism condition for quantum Heisenberg manifolds.}
  \author{Beatriz Abadie \and Ruy Exel\thanks{Partially supported by CNPq,
Brazil.}}
  \date{\today}

\begin{document}
  \maketitle
  \begin{abstract}

 {\bf Abstract:} A study of Hilbert \cstar-bimodules over commutative
$C^*$-algebras is carried out and used to establish a sufficient condition for
two quantum Heisenberg manifolds to be isomorphic.

  \end{abstract} {\bf{\underline{Introduction}.}} In \cite{aee}, a theory of
crossed products of $C^*$-algebras by \hcbs was introduced and used to
describe certain deformations of Heisenberg manifolds constructed by Rieffel
(see \cite{rfhm} and \cite[3.3]{aee}).  This deformation consists of a family
of $C^*$-algebras, denoted $D^c_{\mu\nu}$, depending on two real parameters
$\mu$ and $\nu$, and a positive integer $c$. In case $\mu=\nu=0$, $\dc$ turns
out to be isomorphic to the algebra of continuous functions on the Heisenberg
manifold $M^c$.

 For K-theoretical reasons \cite{fpa}, $D^c_{\mu\nu}$ and $D^{c'}_{\mu'\nu'}$
cannot be isomorphic unless $c=c'$.  It is the main purpose of this work to
show that the \cstar-algebras $D^c_{\mu\nu}$ and $D^c_{\mu'\nu'}$ are
isomorphic when $(\mu,\nu)$ and $(\mu',\nu')$ are in the same orbit under the
usual action of $GL_2(Z)$ on the torus $T^2$ (here the parameters are viewed
as running in $T^2$, since $\dc$ and $D^c_{\mu +n,\nu +m}$ are isomorphic for
any integers $m$ and $n$).

As indicated above, the quantum Heisenberg manifold $D^c_{\mu\nu}$ may be
described as a crossed product of the commutative $C^*$-algebra $\ct$ by a
Hilbert \cstar-bimodule.  Motivated by this, we are led to study some special
features of \hcbs over commutative $C^*$-algebras, which are relevant to our
purposes.

In Section 1 we consider, for a commutative \cstar-algebra $A$, two subgroups
of its Picard group $\pic(A)$: the group of automorphisms of $A$ (embedded in
$\pic(A)$ as in \cite{bgr}), and the classical Picard group $\cpic(A)$ (see,
for instance, \cite{dg}) consisting of Hilbert line bundles over the spectrum
of $A$. Namely, we prove that $\pic(A)$ is the semidirect product of
$\cpic(A)$ by $\aut(A)$. This result carries over a slightly more general
setting, and a similar statement (see Proposition \ref{leftiso}) holds for
\hcbs that are not full, partial automorphisms playing then the role of
$\aut(A)$.  These results provide a tool that enables us to deal with
$\pic(\ct)$ in order to prove our isomorphism theorem for quantum Heisenberg
manifolds, which is done in Section 2.

The authors would like to acknowledge the financial support from FAPESP (grant
no. 95/4609-0), Brazil, and CONICYT, Uruguay.

  \section{The Picard group and the classical Picard group.}
{\bf{\underline{Notation}.}}  Let $A$ be \cstar-algebra. If $M$ is a \hcb over
$A$ (in the sense of \cite[1.8]{bms}) we denote by $\langle$ , $\lrangle ^M$,
and $\langle$ , $\rrangle^M$, respectively, the left and right $A$-valued
inner products, and drop the superscript whenever the context is clear
enough. If $M$ is a left (resp. right) Hilbert \cstar-module over $A$, we
denote by $K(_{A}\!M)$ (resp. $K(M_{A})$) the \cstar-algebra of compact
operators on $M$.  When $M$ is a \hcb over $A$ we will view the elements of
$\langle M, M\rrangle$ (resp. $\langle M, M\lrangle$) as compact operators on
the left (resp. right) module $M$, as well as elements of $A$, via the
well-known identity:
  \[\langle m,n\lrangle p= m \langle n,p \rrangle,\] for $m,n,p\in M$.

The bimodule denoted by $\tilde{M}$ is the dual bimodule of $M$, as defined in
\cite[6.17]{irep}.

By an isomorphism of left (resp. right) Hilbert \cstar-modules we mean an
isomorphism of left (resp. right) modules that preserves the left
(resp. right) inner product. An isomorphism of \hcbs is an isomorphism of both
left and right Hilbert \cstar-modules.  We recall from \cite[3]{bgr} that
$\pic(A)$, the Picard group of $A$, consists of isomorphism classes of full
\hcbs over $A$ (that is, \hcbs $M$ such that $\langle M,M \lrangle =\langle
M,M \rrangle = A$), equipped with the tensor product, as defined in
\cite[5.9]{irep}.

It was shown in \cite[3.1]{bgr} that there is an anti-homomorphism from
$\aut(A)$ to $\pic(A)$ such that the sequence
  \[ 1\longrightarrow \gin(A) \longrightarrow \aut(A) \longrightarrow
\pic(A)\] is exact, where $\gin(A)$ is the group of generalized inner
automorphisms of $A$.  By this correspondence, an automorphism $\alpha$ is
mapped to a bimodule that corresponds to the one we denote by
$A_{\alpha^{-1}}$ (see below), so that $\alpha\mapsto A_{\alpha}$ is a group
homomorphism having $\gin(A)$ as its kernel.

Given a partial automorphism $(I,J,\theta)$ of a \cstar-algebra $A$, we denote
by $J_{\theta}$ the corresponding (\cite[3.2]{aee}) \hcb over $A$.  That is,
$J_{\theta}$ consists of the vector space $J$ endowed with the $A$-actions:
  \[a\cdot x= ax \comsp x\cdot a = \theta[\theta^{-1}(x)a],\] and the inner
products
  \[\langle x,y \lrangle=xy^*,\] and
  \[\langle x,y \rrangle=\theta^{-1}(x^*y),\] for $x,y\in J$, and $a\in A$.
If $M$ is a \hcb over $A$, we denote by $M_{\theta}$ the \hcb obtained by
taking the tensor product $M\otimes_A J_{\theta}$.

The map $m\otimes j \mapsto mj$, for $m\in M,$ $j\in J$, identifies
$M_{\theta}$ with the vector space $MJ$ equipped with the $A$-actions:
  \[a\cdot mj= amj \comsp mj\cdot a =m \theta[\theta^{-1}(j) a],\] and the
inner products
  \[ \langle x,y\lrangle^{M_{\theta}}=\langle x,y\rangle^{M}_L,\] and
  \[\langle x,y\rrangle^{M_{\theta}} =\theta^{-1}(\langle x,y \rangle^M_R),\]
where $m \in M$, $j\in J$, $x,y\in MJ$, and $a\in A$.

As mentioned above, when $M$ is a \cstar-algebra $A$, equipped with its usual
structure of \hcb over $A$, and $\theta\in \aut(A)$ the bimodule $A_{\theta}$
corresponds to the element of $\pic(A)$ denoted by $X_{\theta^{-1}}$ in
\cite[3]{bgr}, so we have $A_{\theta}\otimes A_{\sigma}\cong A_{\theta\sigma}$
and $\widetilde{A_{\theta}}\cong A_{\theta^{-1}}$ for all $\theta,\sigma\in
\aut(A)$.

In this section we discuss the interdependence between the left and the right
structure of a Hilbert \cstar-bimodule. Proposition \ref{leftiso} shows that
the right structure is determined, up to a partial isomorphism, by the left
one.  By specializing this result to the case of full \hcbs over a commutative
\cstar-algebra, we are able to describe $\pic(A)$ as the semidirect product of
the classical Picard group of $A$ by the group of automorphisms of $A$.

  \begin{prop}
  \label{leftiso}
 Let $M$ and $N$ be Hilbert
  \cstar-bimodules over a \cstar-algebra $A$. If $\Phi:M\longrightarrow N$ is
an isomorphism of left $A$-Hilbert \cstar-modules, then there is a partial
automorphism $(I,J,\theta)$ of $A$ such that $\Phi:M_{\theta}\longrightarrow
N$ is an isomorphism of $A-A$ Hilbert \cstar-bimodules.  Namely, $I=\langle
N,N\rrangle$, $J=\langle M,M\rrangle$ and $\theta(\langle \Phi(m_0),
\Phi(m_1)\rrangle)=\langle m_0,m_1\rrangle$.
  \end{prop}

  \pf Let $\Phi: M\longrightarrow N$ be a left $A$-Hilbert \cstar-module
isomorphism.  Notice that, if $m\in M$, and $\|m\|=1$, then, for all $m_i,m'_i
\in M$, and $i=1,...,n$:

  \[\begin{array}[t]{ll}

  \|\sum m\langle m_i,m'_i\rrangle\| & =\|\sum\langle m,m_i\lrangle m'_i\|\\ &
=\|\Phi(\sum\langle m,m_i\lrangle m'_i)\|\\ & =\|\sum\langle m,m_i\lrangle
\Phi(m'_i)\|\\ &=\|\sum \langle \Phi(m),\Phi(m_i)\lrangle \Phi(m'_i)\|\\
&=\|\sum \Phi(m) \langle \Phi(m_i),\Phi(m'_i)\rrangle \|.
  \end{array}
  \] Therefore:
  \[\begin{array}[t]{ll}
  \|\sum \langle m_i,m'_i \rrangle\|&=\sup_{\{m:\|m\|=1\}}\|\sum m\langle
m_i,m'_i\rrangle\|\\ &=\sup_{\{m:\|m\|=1\}}\|\sum \Phi(m) \langle
\Phi(m_i),\Phi(m'_i)\rrangle\|\\ &=\|\sum \langle \Phi(m_i),\Phi(m'_i)
\rrangle\|,
  \end{array}\]

Set $I=\langle N,N\rrangle$, and $J=\langle M,M\rrangle$, and let
$\theta:I\longrightarrow J$ be the isometry defined by
  \[ \theta (\langle \Phi(m_1), \Phi(m_2)\rrangle)= \langle m_1,m_2
\rrangle,\] for $m_1,m_2\in M$.  Then,
  \[\begin{array}[t]{ll}
  \theta(\langle \Phi(m_1), \Phi(m_2) \rrangle^*)&=\theta(\langle \Phi(m_2),
\Phi(m_1) \rrangle)\\ &=\langle m_2, m_1\rrangle\\ &=\langle
m_1,m_2\rrangle^*\\ &=\theta(\langle \Phi(m_1), \Phi(m_2) \rrangle)^*,
  \end{array}\] and
  \[\begin{array}[t]{ll}
  \theta(\langle \Phi(m_1), \Phi(m_2) \rrangle \langle \Phi(m'_1), \Phi(m'_2)
\rrangle)&=\theta(\langle \Phi(m_1), \Phi(m_2)\langle \Phi(m'_1), \Phi(m'_2)
\rrangle \rrangle)\\ &=\theta(\langle \Phi(m_1), \langle \Phi(m_2),
\Phi(m'_1)\lrangle \Phi(m'_2) \rrangle)\\ &=\langle
m_1,\langle\Phi(m_2),\Phi(m'_1)\lrangle m'_2\rrangle\\ &=\langle m_1,\langle
m_2,m'_1\lrangle m'_2\rrangle\\ &=\langle m_1,m_2\langle m'_1,m'_2\rrangle
\rrangle\\ &=\langle m_1,m_2\rrangle \langle m'_1,m'_2\rrangle\\
&=\theta(\langle m_1,m_2\rrangle)\theta(\langle m'_1,m'_2\rrangle) ,
  \end{array}\] which shows that $\theta$ is an isomorphism.

Besides, $\Phi:M_{\theta}\longrightarrow N$ is a \hcb isomorphism:

  \[\begin{array}[t]{ll}
  \Phi(m\langle m_1,m_2\rrangle\cdot a)&=\Phi(m\theta[\theta^{-1}(\langle
m_1,m_2\rrangle)a]\\ &=\Phi(m\theta(\langle \Phi(m_1),\Phi(m_2)a\rrangle))\\
&=\Phi(m\langle m_1,\Phi^{-1}(\Phi(m_2)a)\rrangle)\\ &=\Phi(\langle
m,m_1\lrangle\Phi^{-1}(\Phi(m_2)a))\\ &=\langle m,m_1\lrangle\Phi(m_2)a\\
&=\Phi(\langle m,m_1 \lrangle m_2)a\\ &=\Phi(m\langle m_1,m_2\rrangle)a,
  \end{array}\] and
  \[\langle \Phi(m_1),\Phi(m_2)\rrangle=\theta^{-1}(\langle
m_1,m_2\rrangle^M)=\langle m_1,m_2\rrangle^{M_{\theta}}.\]

Finally, $\Phi$ is onto because
  \[\Phi(M_{\theta})=\Phi(M\langle M,M\rrangle)=\Phi(M)=N.\]

  \qed
  \vspace{.2in}
  \begin{clly}
  \label{unrs} Let $M$ and $N$ be \hcbs over a \cstar-algebra $A$, and let
$\Phi:M\longrightarrow N$ be a an isomorphism of left Hilbert
  \cstar-modules.  Then $\Phi$ is an isomorphism of \hcbs if and only if
$\Phi$ preserves either the right inner product or the right $A$-action.
  \end{clly}
  \pf Let $\theta$ be as in Proposition \ref{leftiso}, so that
$\Phi:M_{\theta}\longrightarrow N$ is a \hcb isomorphism. If $\Phi$ preserves
the right inner product, then $\theta$ is the identity map on $\langle
M,M\rrangle$ and $M_{\theta}=M$.

If $\Phi$ preserves the right action of $A$, then, for $m_0,m_1,m_2\in M$ we
have:
  \[\begin{array}[t]{ll}
  \Phi(m_0)\langle \Phi(m_1),\Phi(m_2)\rrangle & = \langle
\Phi(m_0),\Phi(m_1)\lrangle \Phi(m_2)\\ &=\langle m_0,m_1\lrangle \Phi(m_2)\\
&=\Phi(m_0\langle m_1,m_2\rrangle)\\ &=\Phi(m_0)\langle m_1,m_2\rrangle,
  \end{array}\] so $\Phi$ preserves the right inner product as well.
  \qed

  \begin{prop}
  \label{lmod} Let $M$ and $N$ be left Hilbert \cstar-modules over a
\cstar-algebra $A$. If $M$ and $N$ are isomorphic as left $A$-modules, and $
K(_A\!M)$ is unital, then $M$ and $N$ are isomorphic as left Hilbert
\cstar-modules.
  \end{prop}
  \pf First recall that any $A$-linear map $T:M\longrightarrow N$ is
adjointable. For if $m_i,m'_i\in M$, $i=1,...,n$ are such that $\sum\langle
m_i,m'_i\rrangle=1_{K(_A\!M)}$, then for any $m\in M$:
  \[T(m)=T(\sum\langle m,m_i\lrangle m'_i)=\sum\langle m,m_i\lrangle
T(m'_i)=(\sum\xi_{m_i,Tm'_i})(m),\] wher $\xi_{m,n}:M\longrightarrow N$ is the
compact operator (see, for instance, \cite[1]{la}) defined by
$\xi_{m,n}(m_0)=\langle m_0,m\lrangle n$, for $m\in M$, and $n\in N$, which is
adjointable.  Let $T:M\longrightarrow N$ be an isomorphism of left modules,
and set $S:M\longrightarrow N$, $S=T(T^*T)^{-1/2}$.  Then $S$ is an $A$-linear
map, therefore adjointable. Furthermore, $S$ is a left Hilbert \cstar-module
isomorphism: if $m_0,m_1\in M$, then
  \[\begin{array}[t]{ll}
  \langle S(m_0),S(m_1)\lrangle&=\langle
T(T^*T)^{-1/2}m_0,T(T^*T)^{-1/2}m_1\lrangle\\ &=\langle
m_0,(T^*T)^{-1/2}T^*T(T^*T)^{-1/2}m_1\lrangle\\ &=\langle m_0,m_1\lrangle.
  \end{array}\]

  \qed

We next discuss the Picard group of a \cstar-algebra $A$. Proposition
\ref{leftiso} shows that the left structure of a full \hcb over $A$ is
determined, up to an isomorphism of $A$, by its left structure.

This suggests describing $\pic(A)$ in terms of the subgroup $\aut(A)$ together
with a cross-section of the equivalence classes under left Hilbert
\cstar-modules isomorphisms. When $A$ is commutative there is a natural choice
for this cross-section: the family of symmetric \hcbs (see Definition
\ref{sym}).  That is the reason why we now concentrate on commutative
\cstar-algebras and their symmetric Hilbert \cstar-bimodules.

  \begin{prop}
  \label{lrip} Let $A$ be a commutative \cstar algebra and $M$ a \hcb over
$A$.  Then $\langle m,n\lrangle p=\langle p,n\lrangle m$ for all $m,n,p\in M$.
  \end{prop}
  \pf We first prove the proposition for $m=n$, the statement will then follow
{}from polarization identities.

Let $m,p\in M$, then:
  \[\begin{array}[t]{l}
  \langle \langle m,m\lrangle p-\langle p, m\lrangle m,\ \langle m,m\lrangle
p-\langle p, m\lrangle m \lrangle\\
  \\ =\langle \langle m,m\lrangle p,\langle m,m\lrangle p\lrangle -\langle
\langle m,m\lrangle p, \langle p,m\lrangle m\lrangle\\ -\langle \langle
p,m\lrangle m, \langle m,m\lrangle p\lrangle +\langle \langle p,m\lrangle m,
\langle p,m\lrangle m\lrangle\\
  \\ = \langle m\langle m,p\rrangle\langle p ,m\rrangle,m\lrangle -\langle
m,m\lrangle \langle p,m\lrangle \langle m,p\lrangle\\ -\langle p,m\lrangle
\langle m,p\lrangle\langle m,m\lrangle + \langle p,m\lrangle \langle
m,m\lrangle\langle m,p\lrangle\\
  \\ = \langle m\langle p,m\rrangle\langle m,p\rrangle,m\lrangle -\langle
m,m\lrangle \langle p,m\lrangle \langle m,p\lrangle\\
  \\ =\langle m\langle p,m\rrangle, m\langle p,m\rrangle\lrangle -\langle
m,m\lrangle \langle p,m\lrangle \langle m,p\lrangle\\
  \\ =\langle\langle m,p \lrangle m,\langle m,p \lrangle m\lrangle -\langle
m,m\lrangle \langle p,m\lrangle \langle m,p\lrangle\\
  \\ =0.
  \end{array}\] Now, for $m,n,p \in M$, we have:
  \[\begin{array}[t]{ll}
  \langle m,n\lrangle p&=\frac{1}{4}\sum_{k=0}^3i^k \langle m+i^kn,
m+i^kn\lrangle p\\
 & \\ &=\frac{1}{4}\sum_{k=0}^3\i^k\langle p,m+i^kn\lrangle (m+i^kn)\\ & \\
&=\frac{1}{4}\sum_{k=0}^3 i^kp\langle m+i^kn, m+i^kn\rrangle\\ & \\ &=p\langle
n,m\rrangle\\ & \\ &=\langle p,n \lrangle m.
  \end{array}\]

  \vspace{.2in}
  \begin{defi}
  \label{sym}Let $A$ be a commutative \cstar-algebra. A \hcb $M$ over $A$ is
said to be {\em{symmetric}} if $am=ma$ for all $m\in M$, and \mbox{$a\in A$.}
If $M$ is a \hcb over $A$, the {\em{symmetrization}} of $M$ is the symmetric
\hcb $M^s$, whose underlying vector space is $M$ with its given left
Hilbert-module structure, and right structure defined by:

  \[m\cdot a =am \comsp \langle m_0,m_1\rrangle^{M^s}=\langle
m_1,m_0\lrangle^M,\] for $a\in A$, $m,m_0,m_1\in M^s$.  The commutativity of
$A$ guarantees the compatibility of the left and right actions. As for the
inner products, we have, in view of Proposition \ref{lrip}:

  \[\begin{array}[t]{ll}
  \langle m_0, m_1\lrangle^{M^s}\cdot m_2&=\langle m_0,m_1\lrangle ^M m_2\\
&=\langle m_2,m_1\lrangle^M m_0\\ &=m_0\cdot \langle m_2,m_1 \lrangle^M\\
&=m_0\cdot \langle m_1,m_2\rrangle ^{M^s},
  \end{array}\] for all $m_0, m_1,m_2\in M^s$.

  \end{defi}
  \vspace{.2in}
  \begin{rk}
  \label{unsym} By Corollary \ref{unrs} the bimodule $M^s$ is, up to
isomorphism, the only symmetric
  \hcb that is isomorphic to $M$ as a left Hilbert module.
  \end{rk}

  \begin{rk}
  \label{spm} Let $M$ be a symmetric \hcb over a commutative \cstar-algebra
$A$ such that $K(_{A}\!M)$ is unital. By Remark \ref{unsym} and Proposition
\ref{lmod}, a symmetric \hcb over $A$ is isomorphic to $M$ if and only if it
is isomorphic to $M$ as a left module.

  \end{rk}

  \begin{exam}
  \label{proj} Let $A=C(X)$ be a commutative unital \cstar-algebra, and let
$M$ be a \hcb over $A$ that is, as a left Hilbert \cstar-module, isomorphic to
$A^np$, for some $p\in\proj(M_n(A))$. This implies that $pM_n(A)p\cong
K(_A\!M)$ is isomorphic to a \cstar-subalgebra of $A$ and is, in particular,
commutative.  By viewing $M_n(A)$ as $C(X,M_n(\CC))$ one gets that
$p(x)M_n(\CC)p(x)$ is a commutative \cstar-algebra, hence rank $p(x)\leq 1$
for all $x\in X$.

Conversely, let $A=C(X)$ be as above, and let $p:X\longrightarrow
\proj(M_n(\CC))$ be a continuous map, such that rank $p(x)\leq 1$ for all
$x\in X$. Then $A^np$ is a \hcb over $A$ for its usual left structure, the
right action of $A$ by pointwise multiplication, and right inner product given
by:
  \[ \langle m,r\lrangle = \tau(m^*r),\] for $m,r\in A^np$, $a\in A$, and
where $\tau$ is the usual $A$-valued trace on $M_n(A)$ (that is,
$\tau[(a_{ij})]=\sum a_{ii}$).

To show the compatibility of the inner products, notice that for any
\mbox{$T\in M_n(A)$}, and $x\in X$ we have:
  \[ (pTp)(x)=p(x)T(x)p(x)=[\mbox{trace}(p(x)T(x)p(x))]p(x),\] which implies
that $pTp=\tau(pTp)p.$ Then, for $m,r,s \in M$:
  \[\langle m,r\lrangle s=mpr^*sp=m\tau(pr^*sp)p=m\tau(r^*s)=m\cdot \langle
r,s \rrangle.\]

Besides, $A^np$ is symmetric:
  \[\langle m,r\rrangle=\tau(m^*r)=\sum_{i=1}^{n}m_i^*r_i=\langle
r,m\lrangle,\] for $m=(m_1,m_2,...,m_n)$, $r=(r_1,r_2,...r_n)\in M$.

Therefore, by Remark \ref{spm}, if $p,q\in \proj(M_n(A))$, the \hcbs $A^np$
and $A^nq$ described above are isomorphic if and only if $p$ and $q$ are
Murray-von Neumann equivalent.  Notice that the identity of $K(_{A}\!A^np)$ is
$\tau(p)$, that is, the characteristic function of the set $\{x\in X:
\mbox{rank }p(x)=1\}$. Therefore $A^np$ is full as a right module if and only
if rank $p(x)=1$ for all $x\in X$, which happens in particular when $X$ is
connected, and $p\not =0$.

  \end{exam}

  \begin{prop}
  \label{ppic}

Let $A$ be a commutative \cstar-algebra. For any \hcb $M$ over $A$ there is a
partial automorphism $(\langle M,M \rrangle, \langle M,M\lrangle ,\theta)$ of
$A$ such that the map $i: (M^s)_{\theta}\longrightarrow M$ defined by $i(m)=m$
is an isomorphism of \hcbs.
  \end{prop}
  \pf The map $i:M^s\longrightarrow M$ is a left Hilbert \cstar-modules
isomorphism. The existence of $\theta$, with $I=\langle M,M \rrangle$ and
$J=\langle M^s,M^s\rrangle=\langle M,M\lrangle$, follows from Proposition
\ref{leftiso}.

  \qed

We now turn to the discussion of the group $\pic(A)$ for a commutative
\cstar-algebra $A$. For a full \hcb $M$ over $A$, we denote by $[M]$ its
equivalence class in $\pic(A)$. For a commutative \cstar-algebra $A$, the
group $\gin(A)$ is trivial, so the map $\alpha\mapsto A_{\alpha}$ is
one-to-one. In what follows we identify, via that map, $Aut(A)$ with a
subgroup of $\pic(A)$.

Symmetric full \hcbs over a commutative \cstar-algebra $A=C(X)$ are known to
correspond to line bundles over $X$. The subgroup of $\pic(A)$ consisting of
isomorphism classes of symmetric \hcbs is usually called the classical Picard
group of $A$, and will be denoted by \cpic($A$).  We next specialize the
result above to the case of full bimodules.

  \begin{nota} For $\alpha\in \aut(A)$, and $M$ a \hcb over $A$, we denote by
$\alpha(M)$ the \hcb $\alpha(M)=A_{\alpha}\otimes M \otimes A_{\alpha^{-1}}$.
  \end{nota}
  \begin{rk}
  \label{alfam} The map $a\otimes m\otimes b \mapsto amb$ identifies
$A_{\alpha}\otimes M\otimes A_{\alpha^{-1}}$ with $M$ equipped with the
actions:
  \[a\cdot m =\alpha^{-1}(a)m \comsp m\cdot a =m\alpha^{-1}(a),\] and inner
products
  \[\langle m_0,m_1\lrangle=\alpha(\langle m_0,m_1 \lrangle^M),\] and

  \[\langle m_0,m_1\rrangle=\alpha(\langle m_0,m_1 \rrangle^M),\] for $a\in
A$, and $m, m_0,m_1\in M.$
  \end{rk}

  \begin{thm}
  \label{sdpr} Let $A$ be a commutative \cstar-algebra. Then \cpic$(A)$ is a
normal subgroup of $\pic(A)$ and
  \[ \pic(A)=\cpic(A) \cp \aut(A),\] where the action of $\aut(A)$ is given by
conjugation, that is $\alpha\cdot M=\alpha(M)$.
  \end{thm}
  \pf Given $[M]\in \pic(A)$ write, as in Proposition \ref{ppic}, $M\cong
M^s_{\theta}$, $\theta$ being an isomorphism from $\langle M,M \rrangle =A$
onto $\langle M,M \lrangle =A$.

Therefore $M\cong M^s\otimes A_{\theta}$, where $[M^s]\in \cpic(A)$ and
$\theta \in \aut(A)$. If $[S]\in \cpic(A)$ and $\alpha\in \aut(A)$ are such
that $M\cong S\otimes A_{\alpha}$, then $S$ and $M^s$ are symmetric bimodules,
and they are both isomorphic to $M$ as left Hilbert \cstar-modules. This
implies, by Remark \ref{unsym}, that they are isomorphic. Thus we have:
  \[M^s\otimes A_{\theta}\cong M^s\otimes A_{\alpha} \Rightarrow
A_{\theta}\cong \widetilde{M^s}\otimes M^s \otimes A_{\theta}\cong
\widetilde{M^s}\otimes M^s \otimes A_{\alpha}\cong A_{\alpha},\] which implies
(\cite[3.1]{bgr}) that $\theta\alpha^{-1}\in \gin(A)=\{id\}$, so
$\alpha=\theta$, and the decomposition above is unique.

It only remains to show that $\cpic(A)$ is normal in $\pic(A)$, and it
suffices to prove that $[A_{\alpha}\otimes S\otimes A_{\alpha^{-1}}]\in
\cpic(A)$ for all $[S]\in \cpic(A)$, and $\alpha\in \aut(A)$, which follows
{}from Remark \ref{alfam}.

  \qed
  \begin{nota} If $\alpha\in\aut(A)$, then for any positive integers $k,l$, we
still denote by $\alpha$ the automorphism of $M_{k\times l}(A)$ defined by
$\alpha[(a_{i j})]=(\alpha(a_{ij}))$.

  \end{nota}
  \begin{lemma}
  \label{alfap}
 Let $A$ be a commutative unital \cstar-algebra, and $ p\in \proj(M_n(A))$ be
such that $A^np$ is a symmetric \hcb over $A$, for the structure described in
Example \ref{proj}. If $\alpha\in\aut(A)$, then $\alpha(A^np)\cong
A^n\alpha(p).$
  \end{lemma}
  \pf Set $ J:\alpha(A^np)\longrightarrow A^n\alpha(p)\comsp J(m\otimes
x\otimes r)=m\alpha(xr)$, for $m\in A_{\alpha}$, $r\in A_{\alpha^{-1}}$, and
$x\in A^np$.  Notice that
  \[m\alpha(xr)=m\alpha(xpr)=m\alpha(xr)\alpha(p)\in A^n\alpha(p).\] Besides,
if $a\in A$
  \[\begin{array}[t]{ll} J(m\cdot a\otimes x\otimes r)&=J(m\alpha(a)\otimes
x\otimes r)\\ &=m\alpha(axr)\\ &=J(m\otimes a\cdot x\otimes r),
  \end{array}\] and
  \[\begin{array}[t]{ll} J(m\otimes x\cdot a \otimes r)&=m\alpha(xar)\\
&=J(m\otimes x\otimes a\cdot r),
  \end{array}\] so the definition above makes sense. We now show that $J$ is a
\hcb isomorphism. For $m\in A_{\alpha}$, $n\in A_{\alpha^{-1}}$, $x\in A^np$,
and $a\in A$, we have:
  \[\begin{array}[t]{ll} J(a\cdot(m\otimes x\otimes r))&=J(am\otimes x\otimes
r)\\ &=am\alpha(xr)\\ &=a\cdot J(m\otimes x\otimes r),
  \end{array}\] and
  \[\begin{array}[t]{ll} J(m\otimes x\otimes r\cdot
a)&=m\alpha(x(r\alpha^{-1}(a))\\ &=m\alpha(xr)a\\ &=J((m\otimes x\otimes
r)\cdot a)
  \end{array}\] Finally,
  \[\begin{array}[t]{ll}
  \langle J(m\otimes x\otimes r),J(m'\otimes x'\otimes r')\lrangle&=\langle
m\alpha(xr),m'\alpha(x'r')\lrangle\\ &=\langle m\cdot
[(xr)(x'r')^*],m'\lrangle\\ &=\langle m\cdot \langle x\cdot\langle
r,r''\lrangle^A,x'\lrangle ^{A^np},m'\lrangle\\ &=\langle m\cdot \langle
x\otimes r,x'\otimes r'\lrangle^{A^np\otimes A_{\alpha^{-1}}},m'\lrangle\\
&=\langle m\otimes x\otimes r,m'\otimes x'\otimes r'\lrangle,
  \end{array}\] which shows, by Corollary \ref{unrs}, that $J$ is a \hcb
isomorphism.

  \qed

  \begin{prop}
  \label{alfa} Let $A$ be a commutative unital \cstar-algebra and $M$ a \hcb
over $A$. If $\alpha\in \aut(A)$ is homotopic to the identity, then
  \[A_{\alpha}\otimes M\cong M \otimes A_{\gamma^{-1}\alpha\gamma},\] where
$\gamma\in\aut(A)$ is such that $M\cong (M^s)_{\gamma}$.
  \end{prop}

  \pf We then have that $K(_A\!M)$ is unital so, in view of Proposition
\ref{lmod} we can assume that $M^s=A^np$ with the Hilbert \cstar-bimodule
structure described in Example \ref{proj}, for some positive integer $n$, and
$p\in \proj(M_n(A))$. Since $p$ and $\alpha(p)$ are homotopic, they are
Murray-von Neumann equivalent (\cite[4]{bl}).  Then, by Lemma \ref{alfap} and
Example \ref{proj}, we have
  \[A_{\alpha}\otimes M\cong A_{\alpha}\otimes M^s\otimes A_{\gamma}\cong
M^s\otimes A_{\alpha\gamma}\cong M\otimes A_{\gamma^{-1}\alpha\gamma}.\]
  \qed

We turn now to the discussion of crossed products by Hilbert \cstar-bimodules,
as defined in \cite{aee}. For a \hcb $M$ over a \cstar-algebra $A$, we denote
by $A\cp_{M} \ZZ$ the crossed product \cstar-algebra. We next establish some
general results that will be used later.
  \begin{nota} In what follows, for $A-A$ Hilbert \cstar-bimodules $M$ and $N$
we write $M\cpcong N$ to denote $A\cp_{M} \ZZ \cong A\cp_{N} \ZZ $.
  \end{nota}

  \begin{prop}
  \label{cpiso} Let $A$ be a \cstar-algebra, $M$ an $A-A$ Hilbert
\cstar-bimodule and $\alpha\in\aut (A)$. Then

i) $M\cpcong \tilde{M}.$

ii) $M\cpcong \alpha(M)$.

  \end{prop}
  \pf Let $i_A$ and $i_M$ denote the standard embeddings of $A$ and $M$ in
$A\cp_{M} \ZZ$, respectively.

{\em{i)}} Set
  \[i_{\tilde{M}}:\tilde{M}\longrightarrow A\cp_{M} \ZZ \comsp
i_{\tilde{M}}(\tilde{m})=i_M(m)^*.\]

Then $(i_A,i_{\tilde{M}})$ is covariant for $(A,\tilde{M})$:
  \[i_{\tilde{M}}(a\cdot \tilde{m})=i_{\tilde
{M}}(\widetilde{ma^*})=[i_M(ma^*)]^*=i_A(a)i_M(m)^*=i_A(a)i_{\tilde{M}}(\tilde{m}),\]
  \[i_{\tilde{M}}(\tilde{m_1})i_{\tilde{M}}(\tilde{m_2})^*=i_M(m_1)^*i_M(m_2)=i_A(\langle
m_0,m_1\rrangle^M)=i_A(\langle m_0,m_1\lrangle^{\tilde{M}}),\] for $a\in A$
and $m,m_0,m_1\in M$. Analogous computations prove covariance on the right.
By the universal property of the crossed products there is a homomorphism from
$A\cp_{{\tilde{M}}} \ZZ$ onto $A\cp_{M} \ZZ$. Since $\tilde{\tilde{M}}=M$, by
reversing the construction above one gets the inverse of $J$.

{\em{ii)}} Set
  \[j_A: A\longrightarrow A\cp_{M} \ZZ \comsp j_{\alpha(M)}: M\longrightarrow
A\cp_{M} \ZZ ,\] defined by $j_A=i_A\comp\alpha^{-1}$ ,
$j_{\alpha(M)}(m)=i_M(m)$, where the sets $M$ and $\alpha(M)$ are identified
as in Remark \ref{alfam}.  Then $(j_A,j_{\alpha(M)})$ is covariant for $(A,
\alpha(M))$:
  \[j_{\alpha(M)}(a\cdot
m)=j_{\alpha(M)}(\alpha^{-1}(a)m)=i_A(\alpha^{-1}(a))i_{M}(m)=j_A(a)i_{\alpha(M)}(m),\]
  \[j_{\alpha(M)}(m_0)j_{\alpha(M)}(m_1)^*=i_M(m_0)i_M(m_1)^*=i_A(\langle
m_0,m_1\lrangle ^M)=\]
  \[=j_A(\alpha\langle m_0,m_1\lrangle^M)=j_A(\langle m_0,m_1\lrangle
^{\alpha(M)}),\] for $a\in A$, $m, m_0,m_1 \in M$, and analogously on the
right.  Therefore there is a homomorphism
  \[J:A\cp_{\alpha(M)} \ZZ\longrightarrow A\cp_{M} \ZZ ,\] whose inverse is
obtained by applying the construction above to $\alpha^{-1}$.

  \qed

  \section{An application: isomorphism classes for quantum Heisenberg
manifolds.}

For $\mu,\nu\in\RR$ and a positive integer $c$, the quantum Heisenberg
manifold $\dc$ (\cite{rfhm}) is isomorphic (\cite[Ex.3.3]{aee}) to the crossed
product $C(\TT ^2)\cp_{(X^c_{\nu})_{\al}} \ZZ,$ where $X^c_{\nu}$ is the
vector space of continuous functions on $\RR\times \TT$ satisfying
$f(x+1,y)=e(-c(y-\nu))f(x,y)$. The left and right actions of $\ct$ are defined
by pointwise multiplication, the inner products by $ \langle
f,g\lrangle=f\overline{g}$, and $\langle f,g\rrangle=\overline{f}g$, and
$\al\in \aut (C(\TT ^2))$ is given by $\al(x,y)=(x+2\mu, y+2\nu)$, and, for
$t\in\RR$, $e(t)=exp(2\pi it)$.  {

Our purpose is to find isomorphisms in the family
  \mbox{$\{\dc: \mu,\nu\in \RR, c\in \ZZ, c > 0\}$}.  We concentrate in fixed
values of $c$, because $K_0(\dc)\cong \ZZ^3 \oplus \ZZ_c$(\cite{fpa}).
Besides, since $\al = \alpha_{\mu +m , \nu + n}$ for all $m,n\in \ZZ$, we view
{}from now on the parameters $\mu$ and $\nu$ as running in $\TT$.

Let $\mc$ denote the set of continuous functions on $\RR\times\TT$ satisfying
  \newline $f(x+1,y)=e(-cy)f(x,y).$ Then $\mc$ is a \hcb over $\ct$, for
pointwise action and inner products given by the same formulas as in $X^c$.

The map $f\mapsto \tilde{f}$, where $\tilde{f}(x,y)= f(x,y+\nu)$, is a Hilbert
\cstar-bimodule isomorphism between $(X^c_{\nu})_{\al}$ and
${\scriptsize{C(\TT^2)_{\sigma}\otimes \mc\otimes C(\TT^2)_{\rho}}}$, where
$\sigma(x,y)=(x,y+\nu)$, and $\rho(x,y)=(x+2\mu,y+\nu)$.  In view of
Proposition \ref{cpiso} we have:
  \[\dc\cong C(\TT ^2)\cp_{C(\TT^2)_{\sigma}\otimes M^c\otimes
C(\TT^2)_{\rho}} \ZZ \cong \]
  \[\cong C(\TT ^2)\cp_{(M^c)_{\rho\sigma}}\ZZ \cong C(\TT ^2)\cp_{\mca}\ZZ.\]
As a left module over $\scriptstyle{\ct}$, $M^c$ corresponds to the module
denoted by $X(1,c)$ in \cite[3.7]{rfcan}. It is shown there that $M^c$
represents the element $(1,c)$ of $K_0(\ct)\cong \ZZ^2$, where the last
correspondence is given by $[X]\mapsto (a,b)$, $a$ being the dimension of the
vector bundle corresponding to $X$ and $-b$ its twist. It is also proven in
\cite{rfcan} that any line bundle over $\ct$ corresponds to the left module
$M^c$, for exactly one value of the integer $c$, and that $M^c\otimes M^d$ and
$M^{c+d}$ are isomorphic as left modules.  It follows now, by putting these
results together, that the map $c\mapsto [M^c]$ is a group isomorphism from
$\ZZ$ to $\cpic(\ct)$.

  \begin{lemma}
  \label{cpct}
  \[\pic(C(\TT^2))\cong \ZZ\cp_{\delta} \aut(C(\TT^2)),\] where
$\delta_{\alpha}(c)=\mbox{{\em{det}}}\alpha_*\cdot c,$ for
$\alpha\in\aut(C(\TT^2))$, and $c\in \ZZ$; $\alpha_*$ being the usual
automorphism of $K_0(C(\TT^2))\cong \ZZ^2$, viewed as an element of
$GL_2(\ZZ)$.
  \end{lemma}

  \pf By Theorem \ref{sdpr} we have:
  \[\pic(\ct)\cong \cpic(\ct)\cp_{\delta}\aut(\ct).\] If we identify
$\cpic(\ct)$ with $\ZZ$ as above, it only remains to show that
$\alpha(M^c)\cong M^{\mbox{det}\alpha_*\cdot c}$. Let us view $\alpha_*\in
GL_2(\ZZ)$ as above.  Since $\alpha_*$ preserves the dimension of a bundle,
and takes $\ct$ (that is, the element $(1,0)\in \ZZ^2$) to itself, we have
  \[\alpha_*=\left( \begin{array}{cc}
        1 & 0 \\
        0 & \mbox{det} \alpha_*
       \end{array}
        \right)\]

Now,
  \[\alpha_*(M^c)=\alpha_*(1,c)=(1,\mbox{det}\alpha_*\cdot
c)=M^{\mbox{det}\alpha_*\cdot c}.\] Since there is cancellation in the
positive semigroup of finitely generated projective modules over $\ct$
(\cite{rfcan}), the result above implies that $\alpha_*(M^c)$ and
$M^{\mbox{det}\alpha_*\cdot c}$ are isomorphic as left modules. Therefore, by
Remark \ref{spm}, they are isomorphic as \hcbs.

  \qed

  \vspace{.2in}

  \begin{thm}
  \label{qhmiso}
 If $(\mu,\nu)$ and $(\mu',\nu')$ belong to the same orbit under the usual
action of $GL(2,\ZZ)$ on $\TT^2$, then the quantum Heisenberg manifolds $\dc$
and $ D^{c}_{\mu'\nu'}$ are isomorphic.
  \end{thm}

  \pf If $(\mu,\nu)$ and $(\mu',\nu')$ belong to the same orbit under the
action of $GL(2,\ZZ)$, then $\alpha_{\mu'\nu'}=\sigma\al\sigma^{-1}$, for some
$\sigma\in GL(2,\ZZ)$. Therefore, by Lemma \ref{cpct} and Proposition
\ref{cpiso}:

  \[\mc_{\alpha_{\mu'\nu'}}\cong \mc_{\sigma\al\sigma^{-1}}\cong M^c\otimes
\ct_{\sigma\alpha\sigma^{-1}}\cong \]
  \[\cong\ct_{\sigma}\otimes M^{\mbox{det}\sigma_*^{-1}\cdot c}\otimes
\ct_{\al\sigma^{-1}}\cong\sigma(M_{\alpha_{\mu\nu}}^{\mbox{det}\sigma_*\cdot
c})\cpcong M^{\mbox{det}\sigma_*\cdot c}_{\al}.\] In case $\mbox{det}
\sigma_*=-1$ we have
  \[M^{\mbox{det}\sigma_*\cdot c}_{\al}\cong M^{-c}_{\al}\cpcong
\widetilde{M^{-c}_{\al}}\cong \ct_{\alpha^{-1}_{\mu\nu}}\otimes M^c\cong
(M^c)_{\alpha^{-1}_{\mu\nu}},\] since $\mbox{det}\alpha_*=1$, because $\al$ is
homotopic to the identity.

On the other hand, it was shown in \cite[0.3]{gpots} that
$M^c_{\alpha_{\mu,\nu}^{-1}}\cpcong\mca$.

Thus, in any case, $M^c_{\alpha_{\mu'\nu'}}\cpcong
M^c_{\alpha_{\mu\nu}}$. Therefore
  \[D^c_{\mu'\nu'}\cong \ct\cp_{M^c_{\alpha_{\mu'\nu'}}}\ZZ\cong \ct\cp
_{M^c_{\alpha_{\mu\nu}}}\ZZ\cong D^c_{\mu\nu}.\]

  \qed

  \noindent {\footnotesize\sc BA: Centro de Matem\'aticas, Facultad de
Ciencias, Universidad de la Rep\'ublica, Eduardo Acevedo 1139, C.P 11 200,
Montevideo, Uruguay.  E-mail address: {\tt abadie@@cmat.edu.uy}}

  \noindent {\footnotesize\sc RE: Departamento de Matem\'atica, Universidade
de S\~ao Paulo, Cidade Universit\'aria "Armando de Salles Oliveira". Rua do
Mat\~ao 1010, CEP 05508-900, S\~ao Paulo, Brazil. E-mail address: {\tt
exel@@ime.usp.br}}

  \end{document}